\documentclass[11pt]{amsart}
\usepackage{amsfonts,amssymb,stmaryrd,amscd,amsmath,latexsym,amsbsy,a4wide}

\theoremstyle{plain}

\newtheorem*{thrm}{Theorem}
\newtheorem{prop}{Proposition}[section]

\newtheorem{lem}{Lemma}[section]
\newtheorem{corr}{Corollary}

\newtheorem*{cor'}{Corollary'}

\newtheorem*{sol}{Solution}

\theoremstyle{definition}
\newtheorem*{defi}{Definition}

\theoremstyle{remark}
\newtheorem{rema}{Remark}[section]
\newcommand{\Diff}{{\rm D}}
\newcommand{\im}{{\rm Im}}

\newcommand{\CC}{\mathbb C}
\newcommand{\ZZ}{\mathbb Z}
\newcommand{\solu}[1]{\begin{sol}{\bf (\ref{#1})}}
\newcommand{\dd}{\partial}

\newcommand{\ord}{{\rm deg}}
\title{Deformed Harish-Chandra homomorphism for the cyclic quiver}
\author{ Alexei Oblomkov}
\address{Department of Mathematics, MIT, 77, Massachusetts Ave., Cambridge,
MA 02139, USA.}
\thanks{This work was partially supported by the NSF grant DMS-9988796 and CRDF grant
RM1-2545-MO-03.}
\email{oblomkov@math.mit.edu}
\date{July 27, 2005}
\begin{document}
\begin{abstract}
In the case of cyclic quiver we prove that the deformed
Harish-Chandra map whose existence was conjectured by Etingof and Ginzburg
is well defined. As an application we prove a Kirillov-type formula for the cyclotomic
Bessel function.
\end{abstract}
\maketitle

\section{Introduction}
In this note we prove  the first part of Conjecture 11.22 from the
paper \cite{EG} on the deformed Harish-Chandra isomorphism for the
cyclic quiver. In other words, we prove that the deformed
Harish-Chandra homomorphism is well defined. The constructed
Harish-Chandra homomorphism can be used to study  the
representation theory of the rational Cherednik algebras
associated to the complex reflection  group $\mathfrak{S}_n\ltimes
(\ZZ_\ell)^n$ where $\mathfrak{S}_n$ is the symmetric group
\cite{Ch,BEG,BEG1,E,GG}. In particular, the paper \cite{G}
explains the construction for the shift functor. Also the last
paper proves the second part of the Conjecture 11.2 concerning the
kernel of the deformed Harish-Chandra homomorphism.

The structure of the text is as follows. In  subsections~\ref{defQ},\ref{defH} we define
the main objects: the cyclic quiver $Q$ along with the space
of representations $RQ_n$ of the associated quiver algebra $\CC Q$, and
the rational Cherednik algebra $H_n(k,c)$ together with its spherical subalgebra
$H_n^{sph}(k,c)$. In  subsections~\ref{defDunkl},\ref{defRad} we define
the Dunkl operator embedding $\Theta_{k,c}$ (and $\Theta_{k,c}^{sph}:
=\Theta_{k,c}|_{H_n^{sph}(k,c)}$) and the radial part map $\mathfrak{R}_{k,c}$ along with
its twisted version $\mathfrak{R}_{k,c}^{tw}$. In  subsection~\ref{Thm} the main theorem
is formulated. It states that the images of $\Theta_{k,c}^{sph}$ and
$\mathfrak{R}^{tw}_{k,c}$ are equal.
Section~\ref{proofs} is devoted to the proof of the theorem. In subsection~\ref{defs}
of
Section~\ref{SBessel} we give a definition of the cyclotomic Bessel function
and in the subsection~\ref{IntForm} we prove an integral  formula for this function. The formula
generalizes the well known
interpretation of the classical Bessel function as
the integral over two dimensional sphere
\cite{V}. In subsection~\ref{sphF} we relate our integral formula to the theory
of spherical functions and Kirillov's orbit method.

The idea of the proof is to establish the equality $\im\Theta_{k,c}^{sph}=\im
\mathfrak{R}^{tw}_{k,c}$ 1) after localization by $x_i=0$, $1\le i\le n$ and 2)
after localization by $x_i^\ell=x_j^\ell$, $1\le i<j\le n$. These two statements
imply the theorem. Statement 1) is checked  using the result for $\ell=1$ proved
in \cite{EG} and 2) is proved  using the result for $n=1$ proved by Holland \cite{H}.

\section{The main theorem}
\subsection{Quiver $Q$}\label{defQ}
Let $Q$ be the cyclic quiver with $\ell$ arrows oriented clockwise. We
label the vertices of the quiver by numbers $0,\dots,\ell-1$ in
the  clockwise direction. Let $RQ_n$ be the space of
representations of the associated quiver algebra $\CC Q$ of dimension
$n\delta$ where $\delta=(1,\dots,1)$. In other words, $RQ_n$ is the
vector space $\oplus_{i=0}^{\ell-1}Hom(V_{i+1},V_{i})$, where $V_i$
is the $n$-dimensional vector space assigned to the  vertex $i$, and
$V_\ell=V_0$. Let us denote by $A_{i,i+1}$ the elements of
$Hom(V_{i+1},V_{i})$.

There is a natural action of the group $G_n=\prod_{i=0}^{\ell-1}
GL(V_i)$ on the space $RQ_n$ by  conjugation. If
$g=(g_0,\dots,g_{\ell-1})\in G$, $g_i\in GL(V_i)$ and
$A=(A_{0,1},\dots, A_{\ell-1,0})$ then $g\cdot A=
(g_0A_{0,1}g_1^{-1},
g_1A_{1,2}g_2^{-1},\dots,g_{\ell-1}A_{\ell-1,0}g_0^{-1})$. Because
the element $(c Id_n,\dots,c Id_n)$ acts trivially on $RQ_n$, we
actually have an action of the group $PG_n=G_n/\mathbb{C}^*$ on
$RQ_n$. We use the notation $\mathfrak{pg}_n$ for the Lie algebra of
$PG_n$.

\subsection{The rational Cherednik algebra $H_n(k,c)$}\label{defH}
We denote by $\Gamma$ the cyclic group $\mathbb{Z}_\ell$ and by
$\Gamma_n$ the semidirect product $\mathfrak{S}_n\ltimes \Gamma^n$ with
the symmetric group $\mathfrak{S}_n$.
The group $\Gamma$ is generated by an element $\alpha$ and the group
$\mathfrak{S}_n$ is generated by transpositions $s_{ij}$ exchanging
$i$ and $j$, $i\ne j$. We
denote by $\gamma_i$ the element of $\Gamma^n$ which is equal to
$\gamma$ at the $i$-th place and $1$ at the other places.
Then the conjugation by the element $s_{ij}$ of the subgroup
$\mathfrak{S}_n\subset \Gamma_n$ acts on the normal subgroup
$\Gamma^n$ by the formula:
\begin{equation*}
s_{ij} \alpha_i^s s_{ij}=\alpha_j^s,\quad s_{ij}\alpha_p^s s_{ij}=\alpha_p^s,
\quad p\ne i,j.
\end{equation*}

Let $L$ be a two dimensional vector space and  fix a basis
$x,y$ in $L$.  Denote by $L_n$ the vector space $L^{\oplus n}$, and let
$x_i,y_i$ denote $x, y$ in the $i$-th component of the sum. The group
$\Gamma_n$ acts on $L_n$:
\begin{gather*}
 \alpha_i(x_i)=\epsilon x_i,\quad\alpha_i(x_j)=x_j,\quad
\alpha_i(y_i)=\epsilon^{-1}y_i,\quad\alpha_i(x_j)=x_j,\\
s_{ij}(x_i)=x_j,\quad s_{ij}(y_i)=y_j\\
s_{ij}(x_k)=x_k,\quad s_{ij}(y_k)=y_k,\quad 1\le i, j\ne k\le n,
\end{gather*}
where $\epsilon$ is a (fixed) primitive $\ell$-th root of unity.

Let $TL_n$ be the tensor algebra of $L_n$, that is the free algebra generated by $x_i,y_i$,
$i=1,\dots,n$. There
is  an action of $\Gamma_n$ on the algebra $TL_n$, hence we can
form a smash product $\Gamma_n\ltimes TL_n$.

\begin{defi}\cite{EG}
The algebra $H_n(k,c)$, $k\in\CC, c\in\CC^{\ell-1}$ is the quotient
of the algebra $\Gamma_n\ltimes TL_n$ by the relations:
\begin{gather*}
[x_i,x_j]=0,\quad [y_i,y_j]=0,\quad 1\le i,j\le n,\\
 [y_i,x_i]=1-k\sum_{j\ne
i}\sum_{m=0}^{\ell-1}s_{ij}\alpha_i^m\alpha_j^{-m}+
\sum_{m=1}^{\ell-1}c_m\alpha_i^m,\quad 1\le i\le n,\\
[y_i,x_j]=k\sum_{m=0}^{\ell-1}s_{ij}\epsilon^m\alpha_i^m\alpha_j^{-m},\quad
1\le i\ne j\le n.
\end{gather*}
If ${\rm e}\in \mathbb{C}[\Gamma_n]$ is the idempotent corresponding
to the trivial representation then $H_n^{sph}(k,c)={\rm e}H_n(k,c)
{\rm e}$ is called the spherical subalgebra of $H_n(k,c)$.
\end{defi}
\subsection{Dunkl  operators}\label{defDunkl} From the definition of $H_n(k,c)$
we see that the subalgebra
generated by $x_i$, $i=1,\dots,n$ is  the algebra of
polynomials of $n$ variables. Thus we can think of $x_i$, $i=1,\dots,n$ as
functions on the $n$-dimensional space
$\mathfrak{t}=\mathbb{C}^n$. Let us denote by $\mathfrak{t}^{reg}$
the open subset of $\mathfrak{t}$ given by the conditions $x_i\ne \epsilon^s x_j$ and
$x_i\ne 0$ for $1\le i\ne j\le n$, $s=0,\dots,\ell-1$.

Let us denote by $\Diff (\mathfrak{t}^{reg})$ the algebra of differential operators on
$\mathfrak{t}^{reg}$ and let  $\Gamma_n\ltimes\Diff (\mathfrak{t}^{reg})$ denote the
smash product. As  explained in \cite{DO} there is an
embedding $\Theta_{k,c}$ of the algebra $H_n(k,c)$ into
$\Gamma_n\ltimes\Diff (\mathfrak{t}^{reg})$. The embedding is given by the formulas
\begin{gather*}
x_i\mapsto x_i,\quad \Gamma_n\ni g\mapsto g, \quad y_i\mapsto
\mathcal{D}_i,\\
\mathcal{D}_i=\frac{\partial}{\partial x_i}+k\sum_{j\ne
i}\sum_{m=0}^{\ell-1}\frac{1}{x_i-\epsilon^m
x_j}(s_{ij}\alpha_i^m\alpha_j^{-m}-1)+
\sum_{m=1}^{\ell-1}\frac{c_m}{(\epsilon^m-1)x_i}(\alpha_i^m-1).
\end{gather*}

The map $\Theta_{k,c}$ induces a map from the spherical
subalgebra $H_n^{sph}(k,c)$ to the algebra $\Diff (\mathfrak{t}^{reg})^{\Gamma_n}$ of the
$\Gamma_n$-invariant differential operators on $\mathfrak{t}^{reg}$. We use the letter
$\Theta^{sph}_{k,c}$ for this map.

\subsection{The radial part map}\label{defRad}
Let us define a character $\chi_c$ of $\mathfrak{pg}_n$ by the
formula
$$ \chi_c(x)=\sum_{i=0}^{\ell-1} C_iTr(x_i),$$
where $C_i=\frac{1}{\ell}-\frac{1}{\ell}
\sum_{m=1}^{\ell-1}\epsilon^{mi}c_m$ for $i=1,\dots,\ell-1$ and
$C_{0}=\frac{1-\ell}{\ell}-\frac{1}{\ell}
\sum_{m=1}^{\ell-1}c_m$. We use the same
letter for the one dimensional representation of $PG_n$.

Let  $pr$: $G_n\to PGL(V_0)$ be  a projection of $G_n$ onto the $0$-th
component composed with the map $GL(V_0)\to PGL(V_0)$. Let $k\in\ZZ_+$,
$W_k\subset\CC[y_1,\dots,y_n]$ is the subspace of
the homogeneous polynomials of degree $kn$.  Let $\rho_k:
\mathfrak{sl}(V_0)\to \mathfrak{gl}(W_k)$ be the corresponding Lie algebra map.
We denote by the same letter $\rho_k$ the
representation of  $\mathfrak{pg}_n$ induced by the
projection $pr$.
 We  use the notation $\rho_{k,c}$ for
$\rho_k\otimes\chi_c$.

To define the radial part map we need the space of
$\mathfrak{pg}_n$-equivariant vector valued functions $Fun_{k,c}$.
Let $Fun'$ be the space of the functions on $RQ_n$ of the form
$f=\tilde{f}\prod_{i=0}^{\ell-1}(\det(A_{i,i+1}))^{r_i}$ where $\tilde{f}$ is a
rational function on $RQ_n$ and $r_i=-\sigma+\sum_{j=0}^iC_j$, $i=0,\dots,\ell-1$,
$\sigma:=\sum_{s=0}^{\ell-1}sC_s/\ell$.
The
function $f\in Fun'\otimes W_k$ is from the space $F_{k,c}$ if
and only if $L_g(f)(x)=\frac{d f(x e^{gt})}{dt}|_{t=0}=\rho_{k,c}(g)f(x)$ for all $g\in
\mathfrak{pg}_n, x\in RQ_n$.\footnote{The fact the functions from $Fun_{k,c}$ are multivalued functions on $RQ_n$
does not create the problem for us because  for the computation of the radial part we only
need the restriction of the function to a small neighborhood of the slice
$\mathcal{S}$ defined below.}

If $x\in \CC^n$ is a vector then denote by  $diag(x)$  the diagonal
matrix of  size $n$.
Let us denote by $\Delta$ the map $\mathbb{C}^n\to RQ_n$ which
sends $x$ to the element $(diag(x),\dots,diag(x))$. The image
$\mathcal{S}$ of $\Delta$ is a slice for the action of $PG_n$ on
$RQ_n$. That is, for a generic element $A\in RQ_n$ there exists an
element $g\in PG_n$ such that $g\cdot A=\Delta(x)\in \mathcal{S}$.
Also, it is easy to see that the element $x$ is unique up to the
action of $\Gamma_n$.

The zero weight space $W_k[0]$ is one dimensional and, the
restriction of $f\in Fun_{k,c}$ to $\mathcal{S}$ takes values in
$W_k[0]$. That is, the restriction $Res(f)$ of the function $f$ can be regarded as
scalar function. Moreover,
$f$ is uniquely determined  by  its restriction $Res(f)$, and if $\bar{f}$ is
$\Gamma_n$-equivariant then there exists a function $f\in Fun_{k,c}$ such that
$\bar{f}=Res(f)$. Thus we
can define the radial part map $\mathfrak{R}_{k,c}$:
$\Diff(RQ_n)^{\mathfrak{pg}_n}\to
\Diff(\mathfrak{t}^{reg})$ by the condition:
$$ Res(D(f))=\mathfrak{R}_{k,c}(D)Res(f),$$
for all $f\in Fun_{k,c}$.
Here we used the notation $\Diff(RQ_n)$ for the space of the
differential operators on  $RQ_n$. It is easy to see that in  fact this map lands in
the subspace $\Diff(RQ_n)^{\mathfrak{pg}_n}$  of
$\mathfrak{pg}_n$-invariant differential operators on $RQ_n$. In particular,
for a function $h\in \mathbb{C}[RQ_n]^{\mathfrak{pg}_n}$ we have
$\mathfrak{R}_{k,c}(h)=Res(h).$

\begin{rema} As explained at \cite{EG} one can generalize the definition of
$\mathfrak{R}_{k,c}$ to the case of any $k\in \CC$. Namely, in this one should consider
the representation $\tilde{W}_k=(y_1\dots y_n)^k\CC_{(0)}[y_1^{\pm 1},\dots,y_n^{\pm 1}]$
where
$\CC_{(0)}[y_1^{\pm 1},\dots,y_n^{\pm 1}]$ is the space of the
Laurent polynomials of degree $0$. If $k$ is a positive integer then $W_k$ is a
subrepresentation of $\tilde{W}_k$, so the two settings are equivalent.
\end{rema}
   It is more convenient to use the twisted version $\mathfrak{R}^{tw}$
of the radial part map:
\begin{gather*}
\mathfrak{R}^{tw}_{k,c}(D):=\delta_{k,c}^{-1}\circ\mathfrak{R}(D)\circ\delta_{k,c},\\
\delta_{k,c}(x):=\delta^{k+1}\delta^{\ell\sigma}_{\Gamma},\quad
\delta_{\Gamma}:=\prod_{i=1}^n x_i,\quad \delta=\prod_{1\le i<j\le
n}(x_i^\ell- x_j^\ell).
\end{gather*}

\subsection{Main result}\label{Thm}

\begin{thrm}
 For all values of $k,c$ we have
 $\im \mathfrak{R}_{k,c}^{tw}=
\im \Theta_{k,c}^{sph}$.
\end{thrm}

When $n=1$ the theorem is a particular case of the results of Holland \cite{H}, and when $\ell=1$ the theorem is proved by
Etingof and Ginzburg \cite{EG} who also conjectured the statement of the theorem
for the general $\ell$. The map
$\mathfrak{HC}_{k,c}:=\mathfrak{R}_{k,c}^{tw}\circ\Theta_{k,c}^{-1}$:
$\Diff(RQ_n)^{\mathfrak{pg}_n}\to {\rm e}H_n(k,c){\rm e}$ whose existence follows
from the theorem is called
the deformed Harish-Chandra homomorphism \cite{EG}.

\section{Proofs}\label{proofs}

From the definition of $H_n(k,c)$ we see that we can localize our algebras
$H_n(k,c)$ by inverting a polynomial of $x_i$ which is preserved by the
action of $\Gamma_n$. For example, the polynomials
$\delta,\delta_\Gamma$ have this property.

Let us also introduce the notation $H_n(k)$  for the rational Cherednik
algebra in the case $\ell=1$. This algebra is the
quotient of $\mathfrak{S}_n\ltimes\mathbb{C}\langle
X_1,Y_1,\dots,X_n,Y_n\rangle$ by the relations:
\begin{gather*}
[Y_i,Y_j]=[X_i,X_j]=0,\quad 1\le i,j\le n,\\
[Y_i,X_j]=ks_{ij},\quad 1\le i\ne j\le n,\\
 [Y_i,X_i]=1-k\sum_{j\ne i} s_{ij}.
\end{gather*}

Then we have the following propositions about the localizations.
\begin{prop} If ${\rm e}_{\Gamma}\in\mathbb{C}[\Gamma^n]\subset\mathbb{C}[\Gamma_n]$
is the idempotent corresponding to the trivial representation, and
$H_n(k)_{(X)}$ is the localization by the variables $X_i$ then we
have the following isomorphism
$$ ({e}_{\Gamma}H_n(k,c){e}_{\Gamma})_{\delta_{\Gamma}}\simeq
H_n(k)_{(X)}.$$ The isomorphism is induced by the embedding ${\bold j}$:
$H_n(k)\hookrightarrow
({e}_{\Gamma}H_n(k,c){e}_{\Gamma})_{\delta_{\Gamma}}$ given by:
\begin{equation*}
Y_i\mapsto \ell^{-1}x_i^{1-\ell}y_i,\quad X_i\mapsto x_i^\ell,
\quad 1\le i\le n.
\end{equation*}
\end{prop}
\begin{proof}  Let $T$ is the $n$ dimensional space with the
coordinates $X_i=x_i^\ell$, $i=1,\dots,n$ and $T^{reg}$ is the open subset defined by
the conditions $X_i\ne X_j$, $X_i\ne 0$, $1\le i\ne j\le n$.
We have an obvious isomorphism $\phi$: $
\Diff(T^{reg})\rtimes \mathfrak{S}_n\to {\rm e}_\Gamma\left(
\Diff (\mathfrak{t}^{reg})\rtimes \Gamma_n\right){\rm e}_\Gamma$.
It is enough to check that $\Theta_{k,c}\circ {\bold j}(Y_i)=\phi\circ\Theta_k(Y_i)$,
$i=1,\dots,n$, where $\Theta_k$ is the Dunkl embedding for $H_n(k)$ (see \cite{EG}).
It can be done by a direct computation.
\end{proof}
\begin{prop} We have an isomorphism $H_n(k,c)_\delta\simeq
H_n(0,c)_{\delta}$ induced by the embedding $
{\bold i}:
H_n(0,c)\hookrightarrow H_n(k,c)_{\delta}$ given by:
\begin{equation*}
y_i\mapsto y_i+k\sum_{j\ne
i}\sum_{m=0}^{\ell-1}\frac{1}{x_i-\epsilon^m
x_j}(s_{ij}\alpha_i^m\alpha_j^{-m}-1),\quad x_i\mapsto x_i,\quad
1\le i\le n.
\end{equation*}
\end{prop}
\begin{proof} The proof follows from the  formulas for the Dunkl operators.
\end{proof}
\begin{rema} If $H_1(c)$ is the rational Cherednik algebra for $n=1$ then
we have $$H_n(0,c)=\mathfrak{S}_n\ltimes H_1(c)^{\otimes n}.$$
\end{rema}

 It is shown in the next two lemmas that we can also make compatible localization on
the space $RQ_n$.

Namely, let $\widetilde{RQ_n}\subset RQ_n$ be the open subset
consisting of the points $A\in RQ_n$ with the  property that  the
maps $A_{i,i+1}$, $i=0,\dots,\ell-1$ are invertible. We have a map
$\pi: \widetilde{RQ_n}\to GL_n=GL(V_0)$ which sends $A\in
\widetilde{RQ_n}$ to the product $A_{0,1}A_{1,2}\dots
A_{\ell-1,0}$. Let $\tilde{i}$: $GL(V_0)\to \widetilde{RQ_n}$ be a
section of $\pi$: $\tilde{i}(X)=(X,1,\dots,1)$. We define a
homomorphism $\mathfrak{s}_{k,c}$:
$\Diff(\widetilde{RQ_n})^{\mathfrak{pg}_n}\to
\Diff(GL_n)^{\mathfrak{sl}_n}$ by the condition
$\mathfrak{s}_{k,c}(D)(\tilde{i}^*(f))=\tilde{i}^*(D(f))$ for all
$f\in Fun_{k,c}$, $D\in
\Diff(\widetilde{RQ_n})^{\mathfrak{pg}_n}$.

As we know
$\Diff(GL_n)^{\mathfrak{sl}_n}=\Diff(\mathfrak{gl}_n)_{loc}^{\mathfrak{sl}_n}$
where the subscript 'loc' stands for localization by the determinant.
In particular the construction for the radial part from the
previous section gives us a map $\mathfrak{R}^{tw}_k:
\Diff(\mathfrak{gl}_n)_{loc}^{\mathfrak{sl}_n}\to
\Diff(T^{reg})^{\mathfrak{S}_n}$. Let
us  introduce a map  $\pi_t$: $\mathfrak{t}\to T$ defined by
$\pi_t^*(X_i)=x_i^\ell$, $i=1,\dots,n$.

It turns out that we can put all these maps into a commutative
diagram:
\begin{lem}\label{loc1} The diagram
$$\begin{CD}
\Diff(\widetilde{RQ_n})^{\mathfrak{pg}_n}@>\mathfrak{R}^{tw}_{k,c}>>
\Diff(\mathfrak{t}^{reg})^{\Gamma_n}@<\Theta^{sph}_{k,c}<<
({\rm e}H_n(k,c) {\rm e})_{\delta_\Gamma}\\
@V\mathfrak{s}_{k,c}VV@A\pi^*_tAA@A{\bold j}AA\\
\Diff(\mathfrak{gl}_n)_{loc}^{\mathfrak{pg}_n}@>\mathfrak{R}^{tw}_{k}>>
\Diff(T^{reg})^{\mathfrak{S}_n}@<\Theta^{sph}_k<< ({\rm
e}H_n(k){\rm e})_{(X)}
 \end{CD}
 $$
 is commutative.
\end{lem}
\begin{proof}
The map $\pi_t^*$ acts on the differential operators by the change
of variables $X_i\mapsto x_i^{\ell}$, $i=1,\dots,\ell$, hence from
the description of the map ${\bold j}$ we see that the right half of the
diagram is commutative. That the left half of the diagram commutes follows
from the definition of the radial part map.
\end{proof}

Let us consider another open subset $RQ_n^0\subset RQ_n$
consisting of points $A\in RQ_n^0$ such that the matrix $\pi(A)$
is diagonalizable with distinct eigenvalues and matrices
$Y_i=A_{i,i+1}\dots A_{i-1,i}\in \mathop{End}(V_i)$, $i\ne 0$ are
nondegenerate. Let $\mathcal{T}\subset RQ_n^0$ be the subset of
diagonal matrices (that is, the matrices $A_{i,i+1}$ are diagonal
for $i=0,\dots,\ell-1$) and let us denote by $i$ the embedding of
$\mathcal{T}$ into $RQ_n^0$.

It is elementary to see that $RQ_n^0= PG_n(\mathcal{T})$. Hence
$i$ induces a map $i^*$:
$\Diff(RQ_n^0)^{\mathfrak{pg}_n}\to\Diff^0$
 where $\Diff^0$ is a algebra of differential operators which
 we describe below.
Let  $\mathcal{O}^{K_n}_\mathcal{T}$ be a ring of invariants of
$K_n=\mathfrak{S}_n\ltimes H$, $H=(\mathbb{C}^*)^{n\ell}$.  The
differential operators preserving $\mathcal{O}_\mathcal{T}^{K_n}$
form a subalgebra $\Diff'\subset \Diff(\mathcal{T})$. The elements
of $\Diff'$ also preserve the subspace
$Fun^0_{k,c}:=Fun_{k,c}|_\mathcal{T}$ and the homomorphism
$r:\Diff'\to End_\CC(Fun^0_{k,c})$ is well defined. The algebra
$\Diff^0$ is the image of $r$.

We also have the radial part map $\mathfrak{R}_c$: $\Diff^0\to
\Diff(\mathfrak{t}^{reg})^{\Gamma_n}$. We actually need a twisted
version of this map
$\mathfrak{R}_c^{tw}=\delta^{-\ell\sigma}_\Gamma\circ
\mathfrak{R}_c\circ \delta^{\ell\sigma}_\Gamma$.

\begin{lem}\label{loc2} The diagram
$$
\begin{CD}
\Diff(RQ_n^0)^{\mathfrak{pg}_n}@>\mathfrak{R}^{tw}_{k,c}>>
\Diff(\mathfrak{t}^{reg})^{\Gamma_n}@<\Theta^{sph}_{k,c}<<
({\rm e}H_n(k,c){\rm e})_{\delta}\\
@Vi^*VV @|@A{\bold i}AA\\
\Diff^0@>\mathfrak{R}^{tw}_c>>\Diff(\mathfrak{t}^{reg})^{\Gamma_n}
@<\Theta^{sph}_{0,c}<<({\rm e}H_n(0,c){\rm e})_{\delta}
\end{CD}
$$
is commutative.
\end{lem}
The proof of the lemma is analogous to the proof of the previous lemma.

 We also use notations $\mathfrak{R}^{tw}_c$ and
$\Theta^{sph}_{c}$ for the  radial part and
Dunkl operator maps in the case $n=1$. The results of Holland imply the following
proposition. It is the simplest case of his theorem and we provide a
proof below.
\begin{prop}\label{rk1} If $n=1$ then $\im\mathfrak{R}^{tw}_c=
\im\Theta^{sph}_c$.
\end{prop}
\begin{proof} By the definition of the Dunkl operators, we have
$\Theta^{sph}_c({\rm e}y^\ell{\rm e})=D'$ where
$D'\in\Diff(\mathfrak{t}^{reg})^{\Gamma}$
A simple computation  with  Dunkl operators shows that:
\begin{equation*}
D'=\left(\frac{\dd}{\dd
x}+
\frac{\tilde{C}_0}{x}\right)\left(\frac{\dd}{\dd
x}+\frac{\tilde{C}_0+\tilde{C}_1}{x}\right)\dots\left(\frac{\dd}{\dd x}+\frac{\tilde{C}_0+
\dots
+\tilde{C}_{\ell-1}}{x}\right),
\end{equation*}
where $\tilde{C}_i=\sum_{m=1}^{\ell-1}\epsilon^{mi}c_i$
The differential operator $D'$ acts  on the space $F$ spanned
by the monomials  $x^r$, $r\in\CC$ and $D'$ is determined (up to a scalar) by
the kernel of the action on $F$. It easy to compute this kernel:
 $\ker D'=\langle x^{a_0},\dots,x^{a_{\ell-1}}\rangle$
where $a_i=-\ell\sum_{s=0}^i C_s$.

 The operator
$$D''=\mathfrak{R}_c\left(\frac{\dd}{\dd A_{0,1}}\frac{\dd}{\dd
A_{1,2}}\dots\frac{\dd}{\dd A_{\ell-1,0}}\right),$$ also acts on $F$. The
function $f=A_{0,1}^{r_{0,1}}\dots A_{\ell-1,0}^{r_{\ell-1,0}}$,
 where $r_{i,i+1}=m/\ell-\sigma+\sum_{s=0}^i C_s$, has the properties that
 $f$ is $\chi_c$-equivariant  and $Res f=x^m$. That is, if $m=\ell\sigma+a_i$ then
 the function $f$ does not depend on $A_{i,i+1}$, hence $D''(f)=0$.
 This implies that $\ker D''=\langle x^{b_0},\dots,
 x^{b_{\ell-1}}\rangle$ where $b_i=\ell\sigma+a_i$.

Because both operators are of order $\ell$ we have proved that $D'= x^{-\ell\sigma}
\circ D''\circ x^{\ell\sigma}$.
This implies the inclusion $ \im\Theta_c^{sph}
\subset \im\mathfrak{R}_c^{tw}$
because ${\rm
e}H_1(c){\rm e}$ is generated by $x^\ell$ and $y^\ell$ \cite{LS}.

To prove the opposite inclusion we need the filtrations on $H_1(c)$,
$\Diff(\mathfrak{t}^{reg})$ and $\Diff(RQ_n)$. These filtrations are defined on the generators:
\begin{gather*}
\ord (y)=1, \quad \ord(x)=\ord(\alpha)=0,\\
\ord\left(\frac{\dd}{\dd x}\right)=1,\quad \ord(x)=0,\\
\ord \left(\frac{\dd }{\dd A_{i,i+1}}\right)=1,\quad \ord (A_{i,i+1})=0,\quad i=0,\dots,\ell-1.
\end{gather*}
It is easy to see that  we have
\begin{gather*}
gr H_1(c)=\ZZ_\ell\ltimes\CC[x,y],\quad gr \Diff(\mathfrak{t}^{reg})=\CC[x^{\pm 1},y],\\
gr \Diff (RQ_n)=\CC [A_{0,1},\dots,A_{\ell-1,0},B_{1,0},\dots,B_{0,\ell-1}],
\end{gather*}
where $B_{i+1,i}$ is the image of $\frac{\dd}{\dd A_{i,i+1}}$, $i=0,\dots,\ell-1$.

The maps $\mathfrak{R}_c^{tw}$ and $\Theta_c$ respect the filtrations. Obviously,
the associated graded map $gr\Theta_c^{sph}$ is just the inclusion $\CC[x,y]^{\ZZ_\ell}
\hookrightarrow \CC[x^{\pm 1},y]$. The map $gr\mathfrak{R}_c^{tw}$ is the restriction map:
$f\mapsto f|_{\mathcal{S}}$, where $\mathcal{S}=\{A_{0,1}=\dots=A_{\ell-1,0},
B_{1,0}=\dots=B_{0,\ell-1}\}$. We have
$ \im gr \mathfrak{R}_c^{tw}= \im gr\Theta_c^{sph}$  by the easiest case of the
main theorem from \cite{Gan}.

Let us remark that we have $I_c\subset \ker\mathfrak{R}^{tw}_c$ where $I_c$ is the ideal
 generated by the elements $A_{i,i+1}\frac{\dd}{\dd A_{i,i+1}}-
 A_{i-1,i}\frac{\dd}{\dd A_{i-1,i}}-
C_i$  $i=0,\dots,\ell-1$. Moreover it is easy see that
$gr I_c=\ker gr \mathfrak{R}_c^{tw}$. Hence we have $gr \im\mathfrak{R}^{tw}_c=
\im gr\mathfrak{R}^{tw}_c$
because we have $gr \ker \mathfrak{R}^{tw}_c\supset gr I_c= \ker gr\mathfrak{R}^{tw}_c$
and $gr \ker\mathfrak{R}^{tw}_c\subset \ker gr \mathfrak{R}^{tw}_c$ for the obvious reasons.
Obviously, $gr \im\Theta_c=\im gr\Theta_c$ (because both maps are injective) hence
we get $gr \im\Theta_c=gr \im\mathfrak{R}^{tw}_c$. Together with the inclusion $\im\Theta_c^{sph}
\subset \im\mathfrak{R}_c^{tw}$ this completes the proof.
\end{proof}

\begin{proof}[Proof of the theorem] From the paper \cite{EG} we know that
$\mathfrak{R}^{tw}_k(\Diff(\mathfrak{gl}_n)^{\mathfrak{pg}_n})=
\Theta_k({\rm e}H_n(k){\rm e})$, hence Lemma~\ref{loc1} implies
that
$\mathfrak{R}^{tw}_{k,c}(\Diff(\widetilde{RQ_n})^{\mathfrak{pg}_n})=
\Theta_{k,c}(({\rm e}H_n(k,c){\rm e})_{\delta_\Gamma})$.

On the other hand the map $i^*$ from Lemma~\ref{loc2} can be shown
to be surjective. The proof is the induction by the order of the
differential operators from $\Diff^0$. In particular, the image of
$i^*$ contains the operators $d_m:=\sum_{i=1}^n
\frac{\partial^m}{\partial
(A_{0,1})_{ii}^m}\dots\frac{\partial^m}{\partial
(A_{\ell-1,0})_{ii}^m}$. The same computation as in the
proposition~\ref{rk1} implies that
$\mathfrak{R}_c^{tw}(d_m)=\Theta^{sph}_{0,c}(\sum_{i=1}^n
y_i^{\ell m})$. That  imply
$\mathfrak{R}_{k,c}(\Diff(RQ_n^0))\supset\Theta_{k,c} (({\rm
e}H_n(k,c){\rm e})_{\delta})$ because the algebra ${\rm e}
H_n(k,c){\rm e}$ is generated by $\sum_{i=1}^nx_i^{p\ell}$ and
$\sum_{i=1}^ny_i^{\ell q}$, $p,q\ge 0$.

 The slight modification of
the argument from the proposition~\ref{rk1} proves that we
actually have equality
$\mathfrak{R}_{k,c}(\Diff(RQ_n^0))=\Theta_{k,c} (({\rm
e}H_n(k,c){\rm e})_{\delta})$. But ${\rm e}H_n(k,c){\rm e}$ is
free as $\CC[x_1,\dots,x_n]$-module (under the left
multiplication), hence we have
$$({\rm e}H_n(k,c){\rm e})_{\delta_\Gamma}\cap ({\rm e}H_n(k,c){\rm
e})_{\delta}={\rm e}H_n(k,c){\rm e}.$$
Thus we proved the theorem.
\end{proof}

\section{Cyclotomic Bessel function}\label{SBessel}
\subsection{Definitions}\label{defs}
For $P\in\CC[y_1,\dots,y_n]^{\Gamma_n}$ let  us define a
$\Gamma_n$-invariant differential operator  $\mathcal{D}_P:=
\Theta_{k,c}^{sph}(P)=P(\mathcal{D}_1,\dots,\mathcal{D}_n)$. The differential operators
$\mathcal{D}_P$, $P\in\CC[y_1,\dots,y_n]^{\Gamma_n}$  mutually commute and we can
study their common eigenfunctions. The problem makes sense if we replace $\Gamma_n$
by any complex reflection group $W$ and $D_i$ by the corresponding
Dunkl operators \cite{DO}. In the case when $W$ is a Coxeter group
the problem was studied by \cite{Op}. Most of the proofs from \cite{Op} are valid
in the case when $W$ is a complex reflexion group. Below we refer to these proofs.

Let us define some simply connected domain $\mathcal{C}$ inside $\mathfrak{t}^{reg}$.
For that we
 define the set of cuts of $\mathfrak{t}^{reg}$:
\begin{align*}
{\rm cut}_i&=\{x\in\mathfrak{t}^{reg}|\Re x_i=0, \Im x_i>0\},\\
{\rm cut}_{i,j;m}&=\{x\in\mathfrak{t}^{reg}|\Re (x_i/x_j-e^{2\pi\sqrt{-1}m/\ell})=0,
\Im (x_i/x_j-e^{2\pi\sqrt{-1}m/\ell})>0\},\quad 0\le m <\ell',\\
{\rm cut}_{i,j;m}&=\{x\in\mathfrak{t}^{reg}|\Re (x_i/x_j-e^{2\pi\sqrt{-1}m/\ell})=0,
\Im (x_i/x_j-e^{2\pi\sqrt{-1}m/\ell})<0\},\quad  \ell'\le m <\ell,
\end{align*}
where $1\le i<j\le n$ and  $\ell'=\ell/2$ if $\ell$ is even,
 $\ell'=(\ell-1)/2$
 if $\ell$ is odd. Let  ${\rm cuts}:=\cup_{1\le i<j\le n}\cup_{0\le m <\ell}
 {\rm cut}_{ij;m}\cup\cup_{1\le i\le n}{\rm cut}_i$ and
 $\mathcal{C}=\mathfrak{t}^{reg}\setminus\rm{cuts}$.

Let us fix $\lambda\in\CC^n$. It is easy to see that the space $V^{k,c}_\lambda$
of solutions of the system of equations in the domain $\mathcal{C}$:
\begin{equation*}
\mathcal{D}_P f=P(\lambda)f,\quad \forall P\in\CC[y_1,\dots,y_n]^{\Gamma_n},
\end{equation*}
has finite dimension (it is actually equal to $|\Gamma_n|$ for generic $\lambda$ (see
\cite{Op}, Corollary~3.7)).

A function $f\in V_\lambda^{k,c}$ is analytic in $\mathcal{C}$ and can be analytically
continued to $\mathfrak{t}^{reg}$
but the continued function is  multivalued because $\mathfrak{t}^{reg}$ is not
simply connected. Moreover, we can continue $f\in V_\lambda^{k,c}$  to $\mathfrak{t}$
if we
allow the singularities. Let us denote this continuation by $\tilde{f}$.

 Let us assume that $k,C_i\in\mathbb{R}$, $0\le i\le \ell-1$
 (the relation between
$C$ and $c$ was explained at subsection~\ref{defRad}). To simplify the exposition
 we also assume that $k>0$, $\sum_{i=0}^s C_i\ge 0$, $0\le s<  \ell$.
 The general case
 can be treated similarly.

 The functions $\tilde{f}$, $f\in V_\lambda^{k,c}$ could have the singularities.
 Let us assume that $C_i\notin\ZZ$, $i=1,\dots,\ell-1$. Then
 the local analysis shows (see section~7 of \cite{Op}) that for any $j$, $1\le j\le n$ and $f\in V_\lambda^{k,c}$ we
 can present $\tilde{f}$ in the form
 $$\tilde{f}(x)=\sum_{s\in S} x_j^{a_s}F_s(x),$$
 where $S$ is a subset of $\{ 0,\dots,\ell-1\}$, $a_s$ are defined at the proof of
 the proposition~\ref{rk1} and $F_s$, $s\in S$ is a nonzero function analytic at the generic
 point of the divisor $\{x_j=0\}$.  In the case when  the assumption on $C$ does not hold there is a similar presentation
 for $\tilde{f}$ (see section~7 of \cite{Op}) which involves the logothimic functions.

Let us assume that $k\notin\ZZ+\frac12$.
 Then
 it is also possible to show (see section~7 of \cite{Op}) that for any $i,j,m$,
 $1\le i\ne j \le n$, $0\le m<\ell$ and $f\in V_\lambda$ we can present $\tilde{f}$ in
 the form
 $$ \tilde{f}(x)=\sum_{s\in S}(x_i-\epsilon^m x_j)^{2s} F_s(x)+\sum_{t\in T}
 (x_i-\epsilon^m x_j)^{2k+1+2t} G_t(x),$$
 where
 $S$, $T$ are subsets of $\{ 0,1,\dots,\ell'\}$  and $F_s$, $G_t$, $s\in S$, $t\in T$ are non zero functions analytic at
 the generic point of the divisor $\{x_i-\epsilon^m x_j=0\}$. Again if $k\in\ZZ+\frac12$
 then there is a similar expression for $\tilde{f}$ which involves
 the logorithmic functions.

Let us also remark that the real part $\mathfrak{t}_{\mathbb{R}}^{reg}$
of $\mathfrak{t}^{reg}$ is a subset of $\mathcal{C}$. Hence
the restriction $f|_{\mathfrak{t}_\mathbb{R}^{reg}}$
is a well defined single valued function. Let us also denote by $\Lambda^{reg}$
the subset of $\CC^n$ such that $\lambda\in\Lambda^{reg}$ if and only if
$\lambda_i\ne 0$, $\lambda_i^\ell\ne\lambda_j^\ell$, $1\le i\ne j\le n$.

\begin{defi} The cyclotomic Bessel function $B_\lambda^{k,c}$, $\lambda\in\Lambda^{reg}$
is a function from $V_\lambda^{k,c}$ such that
\begin{enumerate}
\item $\tilde{B}_\lambda^{k,c}$ has no singularities (i.e. takes only finite values)   on
$\mathfrak{t}\simeq \CC^n$,
\item $B_\lambda^{k,c}|_{\mathfrak{t}_\mathbb{R}^{reg}}$ is
$S_n$-invariant,
\item $\tilde{B}_\lambda^{k,c}(0)=1$.
\end{enumerate}
\end{defi}

According to propositions 5.6, 6.8 and corollary 7.16
of \cite{Op} the conditions (1), (2), (3)
  define the function
$B_\lambda^{k,c}$  uniquelly.  Moreover, from the results of \cite{Op} it
follows that the function
$\tilde{B}_\lambda^{k,c}$ is a singlevalued function.

\begin{rema} From the assumptions on the parameters $C_i$, $0\le i<\ell$
we see that $a_i< a_{\ell-1}=0$, $0\le i<\ell-1$.
Hence the generic function from $V_\lambda^{k,c}$ has  singularities.
\end{rema}

\begin{rema} When $n=1$ and $\ell=2$ the function $B_\lambda^{k,c}$ is related to
the classical
Bessel function $J_{C_1}$ by the formula
$(\lambda x)^{C_1}B^{k,c}_\lambda(x)/\Gamma(C_1)=J_{C_1}(2x\lambda)$.
\end{rema}

Using the theory of the deformed Harish-Chandra homomorphism we find a Kirillov-type
\cite{K}
integral formula for the  Bessel function. We use  the ideology of the paper
\cite{EFK} to do that.

Let  $\bold{K}_n=U(n)^{\times\ell}\subset G_n$
be the maximal compact subgroup of $G_n$.
Let $Q^{op}$ be a cyclic quiver with $\ell$ vertices
and the edges  oriented counterclockwise.
Let $RQ^{op}_n$ be the space of representations of $\CC Q^{op}$ of dimension
$n\delta$. Let us denote by $\Delta^{op}$ the map which sends the element
$y\in \CC^n$ to the element $(diag(y),\dots,diag(y))\in R Q_n^{op}$.
Let us denote by $\mathcal{O}_\lambda\subset RQ^{op}_n$ the $\bold{K}_n$-orbit
of the element $\Delta^{op}(\lambda)$. The invariant measure on $\bold{K}_n$ induces
a $\bold{K}_n$-invariant measure $d\mu_{\lambda}$ on $\mathcal{O}_\lambda$.
Then the space $V_\lambda:=
L^2(\mathcal{O}_\lambda)$ has a natural structure
of $\bold{K}_n$-module:
\begin{equation*}
 (g\cdot f)(B)=f(g^{-1}\cdot B).
\end{equation*}
The action of $\bold{K}_n$ respects the Hermitian product $\langle f,g\rangle_\lambda:
=\int_{\mathcal{O}_\lambda} f(A)\bar{g}(A) d\mu_{\lambda}(A)$ on $V_\lambda$.

\subsection{}
Now let us assume that $C\in\ZZ^\ell$   and $k\in\ZZ_+$.
Then the map $\rho_{k,c}$ gives us the representation of $\bold{K}_n$ which
we denote by  $W_{k,c}$ (it is isomorphic to $W_k$ as a vector space).
Let $W_{k,c}^\vee$ be its dual.
The vector space $V_\lambda$ is a unitary representation
of $\bold{K}_n$ and we have:

\begin{prop}\label{intertw} If $\lambda\in\Lambda^{reg}$  then there
exists  a unique up to scaling   injective map of $\bold{K}_n$-representations
$\eta_{k,c}^{\lambda}$: $W_{k,c}^\vee\to V_{\lambda}$.
\end{prop}
\begin{proof}
Because of the conditions on $\lambda$ we have $\mathcal{O}_\lambda\simeq
\bold{K}_n/\bold{T}_n$ where $\bold{T}_n$ is a  torus which is equal to
the stabilizer of $\Delta^{op}(\lambda)$. On  the other hand the space  $L^2(\bold{K}_n)$
has the left and right actions of $\bold{K}_n$ and by  the Peter-Weyl
theorem it decomposes
into the direct sum of subrepresentations
$$L^2(\bold{K}_n)=\hat{\oplus}_{V\in \hat{\bold{K}}_n}V^\vee\otimes V,$$ where
$\hat{\bold{K}}_n$ is a notation for the space of all irreducible representations
of $\bold{K}_n$ and $\hat{\oplus}$ is a notation for the completed direct sum.
Hence we have
$$L^2(\mathcal{O}_\lambda)=\hat{\oplus}_{V\in \widehat{\bold{K}}_n}V^\vee\otimes
(V)^{\bold{T}_n}.$$
To finish the proof let us notice that
$W_{k,c}^{\bold{T}_n}\simeq \CC$.
\end{proof}

\subsection{Integral formula}\label{IntForm}
From the explicit construction of $W_{k,c}\subset\CC[y_1,\dots,y_n]$ we know that the
monomials $y_1^{i_1}\dots y_n^{i_n}$, $\sum_{j=1}^n i_j=nk$ span $W_{k,c}$
and $W_{k,c}[0]=\langle (y_1\dots y_n)^k\rangle$.
Let us define the function $m_k$ on $U(n)$ by the formula:
$$ g((y_1\dots y_n)^k)=m_k(g)(y_1\dots y_n)^k+\mbox{the linear combination of the
other monomials}.$$ Then we have
\begin{corr}\label{maincor} For $\lambda\in \Lambda^{reg}$ we have
\begin{equation}\label{BesselInt}
B_{\lambda}^{k,c}(x)\sim\frac{
\int_{\bold{K}_n} e^{tr(g\Delta^{op}(\lambda)g^{-1}\Delta(x))}
m_k(g_0)\prod_{i=0}^{\ell-1}det(g_i)^{C_i} d\mu_l(g)}{\delta_{k,c}(x)},
\end{equation}
where $d\mu_l$ is a left invariant measure on $\bold{K}_n$ and  $\sim$ stands for being
proportional.
\end{corr}
Before giving the prove let us discuss the formula (\ref{BesselInt}) in the case
$n=1$.

\subsection{Case n=1} In this case we can omit $k$ from the notations. The space
$V_\lambda^c$ is a space of solutions of the ODE:
$$\left(\frac{\dd}{\dd
x}+
\frac{\tilde{C}_0}{x}\right)\left(\frac{\dd}{\dd
x}+\frac{\tilde{C}_0+\tilde{C}_1}{x}\right)\dots\left(\frac{\dd}{\dd x}+\frac{\tilde{C}_0+
\dots
+\tilde{C}_{\ell-1}}{x}\right)f=\lambda^l f.$$

Because of the assumption on the parameters $C_i$, $0\le i<\ell$ we know that
$a_i<a_{\ell-1}=0$, $0\le i<\ell$.
Hence
the space of the function from $V_\lambda^c$ without a pole at $x=0$ is one-dimensional and spanned by
the Bessel function $B_\lambda^c$.

Corollary gives the formula:
\begin{equation}\label{rk1IntBessel} B_\lambda^c(x)\sim x^{-\sum_{s=0}^{\ell-1} sC_s}
\int_{[0,1)^\ell}e^{\sum_{j=0}^{\ell-1}
2\pi\sqrt{-1}\varphi_jC_j+e^{2\pi\sqrt{-1}(\varphi_j-\varphi_{j+1})}\lambda x}
d \varphi_0\dots d\varphi_{\ell-1},
\end{equation}
where we assume that $\varphi_\ell=\varphi_0$.
\begin{proof}[Proof of Corollary for $n=1$]
It is easy to check that RHS of   (\ref{rk1IntBessel})
is a function from $V_\lambda^c$.
Let us explain why it has no pole at $x=0$. That is we need to prove that the integral
in   (\ref{rk1IntBessel}) has the zero of order $\sum_{s=0}^{\ell-1} sC_s$ at $x=0$.

Let us denote the integrand in (\ref{rk1IntBessel})  by $F(\varphi,x)$.
Let $e_i$, $i=0,\dots,\ell-1$
be a standard basis in $\ZZ^\ell$ and $\nu_i=e_i-e_{i+1}$, $0\le i<\ell-1$,
$\nu_{\ell-1}=e_{\ell-1}-e_0$. Then  we have:
$$
\frac{\partial^t}{\partial x^t}F(0,\varphi)=\lambda^t
\sum_{i_1,\dots,i_t=0}^{\ell-1}e^{(C-\sum_{s=1}^t\nu_{i_s},\varphi)},$$
where $(y,\varphi)=2\pi\sqrt{-1}\sum_{i=0}^{\ell-1}y_i\varphi_i$.  Only terms such that
$C-\sum_{s=1}^t\nu_{i_s}=0$ give a nonzero input into the integral
$\int_{[0,1)^\ell}\frac{\partial^t}{\partial x^t}F(0,\varphi)d\varphi_0\dots
d\varphi_{\ell-1}$. There are no such terms
if $t< \sum_{s=0}^{\ell-1} sC_s$. Moreover, if $t=\sum_{s=0}^{\ell-1}sC_s$ then there is only one
such term.
\end{proof}

\begin{corr} For $\lambda\ne 0$ we have
\begin{equation*} B_\lambda^c(x)= \frac{t!}{(\lambda x)^t}
\int_{[0,1)^\ell}e^{\sum_{j=0}^{\ell-1}
2\pi\sqrt{-1}\varphi_jC_j+e^{2\pi\sqrt{-1}(\varphi_j-\varphi_{j+1})}\lambda x}
d \varphi_0\dots d\varphi_{\ell-1},
\end{equation*}
where $t=\sum_{s=0}^{\ell-1} sC_s$.
\end{corr}

\subsection{}\begin{proof}[Proof of the Corollary~\ref{maincor}]
One can show that  we actually have $\im \eta_{k,c}^{\lambda}\subset C^{\infty}
(\mathcal{O}_\lambda)$.
Let us choose a nonzero function $\varphi^{k,c}_\lambda\in\eta_{k,c}^{\lambda}(
W_{k}^\vee[0])\in C^{\infty}(\mathcal{O}_\lambda)$ where the space $W_k[0]^\vee\subset W_{k,c}$
is the one defined in subsection~\ref{IntForm}.
Then it is  easy to see that RHS of (\ref{BesselInt}) is proportional to the function:
$$
h_\lambda(x)=\frac{\int_{\mathcal{O}_\lambda}
\varphi^{k,c}_\lambda(A)e^{tr(A\Delta(x))} d\mu_\lambda(A)}{\delta_{k,c}(x)}.
$$

The space $W^\vee_{k,c}$ is a subspace of $V_\lambda$  and we can choose a
basis $v^\vee_1,\dots,v^\vee_N$ in $W^\vee_{k,c}$ and
$v_1,\dots,v_N$ be its dual basis. Let us define the vector valued function
$$f_{\lambda}(B)=\sum_{i=1}^N v_i\int_{\mathcal{O}_\lambda} v^\vee_i(A)exp(tr(AB))
d\mu_\lambda(A).$$
It is easy to see that $f_\lambda\in Fun_{k,c}$
and $f_\lambda|_{\mathfrak{t}}$ is proportional to $h_\lambda$.
For every $P\in\CC[y_1,\dots,y_n]^{\Gamma_n}$
there exists $G_n$-in\-va\-ri\-a\-nt function $S_P\in\CC[R Q_n^{op}]^{PG_n}$ such that
$S_P(\Delta^{op}(y))=P(y)$. One should think of $S_P$ as a differential
operator on $RQ_n$ with  constant coefficients. From the proof of the main theorem one
can see
that  $\mathfrak{R}^{tw}_{k,c}(S_P)=\mathcal{D}_P$.
Also from the formula
for $f_{\lambda}$ one can see that $S_P f_\lambda=
S_P(\Delta^{op}(\lambda))f_\lambda=P(\lambda) f_\lambda$. Finally we use that
$f_\lambda(\Delta(x))\in W_{k}[0]$ for all $x\in \CC^n$.

Let us prove that RHS of (\ref{BesselInt}) is analytic. From the discussion
before the definition of the Bessel function we see that the functions from
$V_\lambda^{k,c}$ can have the poles only along the divisors $x_j=0$, $j=1,\dots,n$.
Hence we only need to show that the order of vanishing of the integral in  RHS of
(\ref{BesselInt}) is at least $\sum_{s=0}^{\ell-1} sC_s$. Let us show it for $j=1$.

Indeed, every element $g\in{\bold K}_n$ can be uniquely presented in the form $g=\tilde{g}A(\varphi)$, where
$\tilde{g}\in SU(n)^{\times \ell}$ and $\varphi\in [0,1)^{\ell}$,
$A_i=diag(e^{2\pi\sqrt{-1}\varphi_i},1,\dots,1)$, $0\le i<\ell$. For  $\varphi\in
[0,1)^\ell$, $x\in \CC^n$ let $\Delta(x,\varphi)$ be the element of $RQ_n$ such that
$\Delta(x,\varphi)_{i,i+1}=diag(e^{2\pi\sqrt{-1}(\varphi_i-\varphi_{i+1})}x_1,x_2,
\dots,x_n)$, $0\le i<\ell$.
We can rewrite the integral
from (\ref{BesselInt}) in the form:
$$Int=\int_{\varphi\in [0,1)^\ell}e^{(\varphi,C)}
\left(\int_{\tilde{G}_n}e^{tr(\tilde{g}\Delta^{op}(\lambda)\tilde{g}^{-1}\Delta(x,\varphi))}
m_k(\tilde{g}_0) d\tilde{\mu}_l(\tilde{g})\right)d\varphi_0\dots d\varphi_{\ell-1},$$
where $\tilde{G}_n=SU(n)^{\times\ell}$ and $d\tilde{\mu}_l$ is the
corresponding left $\tilde{G}_n$-invariant measure. Now let us notice
that when $x_1=0$ the integral:
$$
\int_{\tilde{G}_n}e^{tr(\tilde{g}\Delta^{op}(\lambda)\tilde{g}^{-1}\Delta(x,\varphi))}
m_k(\tilde{g}) d\tilde{\mu}_l(\tilde{g})
$$
does not depend of $\varphi$. Hence the computation of the order of vanishing of $Int$
at $x_1=0$ is essentially the same as in the case $n=1$.
\end{proof}

\begin{rema} It is easy to see that the integral in the formula (\ref{BesselInt})
can be reduced to the integral over $\bold{K}_n/\bold{T}_n$ because the integrand
is $\bold{T}_n$-invariant.
\end{rema}

\subsection{An interpretation in term of spherical functions}\label{sphF}
Let us illuminate the connection with  the theory of
spherical functions and Kirillov's orbit method.
Let $RQ_n(\mathbb{R})\subset RQ_n$ ($RQ^{opp}_n(\mathbb{R})
\subset RQ^{opp}_n$) be the subspace of the real points of $RQ_n$ ($RQ_n^{opp}$).
Let $G_n(\mathbb{R})$ be the real part of $G_n$. As $G_n(\mathbb{R})$
acts on the space $RQ_n(\mathbb{R})$,  we have the group $\bold{H}_n(\mathbb{R})
=RQ_n(\mathbb{R})\rtimes G_n(\mathbb{R})$ with the relations:
$$ (A,g)\cdot (B,h)=(A+g\cdot B,gh),$$
where $A,B\in RQ_n(\mathbb{R})$ and $g,h\in G_n(\mathbb{R})$.

The group $\bold{H}_n(\mathbb{R})$ has a left and right action of the subgroup
$G_n(\mathbb{R})$.
In particular, the space $Fun_{k,c}$ is the space of
$W_{k,c}$-valued functions on $\bold{H}_n(\mathbb{R})$
which are invariant with respect to the left action of $G_n(\mathbb{R})$
and $G_n(\mathbb{R})$-equivariant with respect to the right action.
Hence from the ideology of the paper \cite{EFK} we know that we can
construct the functions from $Fun_{k,c}$ starting from the representations
of $\bold{H}_n(\mathbb{R})$.

Now we  construct unitary representations of
$\bold{H}_n(\mathbb{R})$.
Let $\mathcal{O}^{\mathbb{R}}_\lambda$, $\lambda\in\mathbb{R}^n$,
$\lambda_i\ne \lambda_j$, $1\le i<j\le n$ be the $G_n(\mathbb{R})$-orbit
of $\Delta^{op}(\lambda)$ inside $RQ_n(\mathbb{R})$. Let $d\mu_\lambda$ be a left
$G_n(\mathbb{R})$-invariant measure on $\mathcal{O}^{\mathbb{R}}_\lambda$. Then
the space $V^{\mathbb{R}}_\lambda:=L^2(\mathcal{O}_\lambda^\mathbb{R})$ has natural Hermitian product:
$\langle f,g\rangle_\lambda=\int_{\mathcal{O}_\lambda} f(A)\bar{g}(A)d\mu_\lambda(A),$
and the group $\bold{H}_n(\mathbb{R})$ acts by unitary operators on
$V_\lambda^{\mathbb{R}}$:
$$
((0,g)\cdot f)(A)=f(g^{-1}\cdot A),\quad ((B,1)\cdot f)(A)=exp(2\pi\sqrt{-1}tr(AB))
f(A),
$$
where $g\in G_n(\mathbb{R})$, $A\in RQ_n(\mathbb{R})$.

Let $\tilde{V}_\lambda^{\mathbb{R}}$ be the space of smooth functions on
$\mathcal{O}_\lambda^{\mathbb{R}}$.
The complexifications of the groups $G_n(\mathbb{R})$ and $\bold{K}_n$ coincide. Hence
Proposition~\ref{intertw} implies that there exists a unique up to scaling map of
$G_n(\mathbb{R})$-representations $\eta_{k,c}^\lambda$: $W^\vee_{k,c}\to
\tilde{V}_\lambda^{\mathbb{R}}$, where $W^\vee_{k,c}$ is a representation of
$G_n(\mathbb{R})$ dual to $W_{k,c}$.

Clearly the image of $\eta_{k,c}^\lambda\subset
\tilde{V}_\lambda^{\mathbb{R}}$ does not
lie inside $V_\lambda^{\mathbb{R}}$ but we  can use pairing
$\langle\cdot,\cdot\rangle_\lambda$ for the elements of $\tilde{V}_\lambda^{\mathbb{R}}$ if we
think of the integration over $\mathcal{O}_\lambda$ in the sense of generalized
functions. Keeping it in mind we can write the formula for the spherical function:
\begin{equation}\label{BesselReal} f_\lambda^{\mathbb R}(X)=\sum_{i=1} v_i
\langle \eta_{k,c}^\lambda(v^\vee_i), X\cdot
1\rangle_\lambda,
\end{equation}
where $X\in \bold{H}_n(\mathbb{R})$, $v_1,\dots, v_N$ is a basis in $W_{k,c}$ and $v^{\vee}_1,\dots,v^\vee_n$ is
a dual basis. It is clear that we have
\begin{corr}For $\lambda\in \mathbb{R}^n$ such that $\lambda_i\ne \lambda_j$,
$1\le i\ne j\le n$ we have
$$f_\lambda^{\mathbb R}((\Delta^{opp}(x),1))\sim
B_\lambda(2\pi\sqrt{-1} x),$$ where $\sim$ stands for being proportional.
\end{corr}
When $\ell=1$ then the coadjoint orbit $\mathcal{O}_\lambda^{\mathbb R}$
corresponds to the  principal series representation $\mathcal{P}^{\lambda}$ of
$GL_n(\mathbb{R})$ (see the User's Guide in \cite{K}). For $k=0$ the function
$f_\lambda^{\mathbb R}((1,A))$,
$A\in\mathfrak{gl}_n(\mathbb{R})$ is the generalized
$\mathfrak{gl}_n(\mathbb{R})$-character of $\mathcal{P}^{\lambda}$ and the formula
(\ref{BesselReal}) is the classical
Kirillov's integral formula for the character. For the general $k$, the construction
is a degeneration of the construction of \cite{EFK} where the vector valued characters
of
$GL(n,\CC)$
were studied. Indeed, the $W_k$-valued character of $G=GL(n,\CC)$ can be
interpreted as function on $G\times G$ which is left $G$-invariant and right
$G$-equivariant (for more detailed discussion see \cite{Op}). While we degenerate
the group $G\times G$ into $\mathfrak{gl}(n,\mathbb{C})\rtimes GL(n,\CC)$ we see that
the construction for the Jack polynomials from \cite{EFK} give us the construction
for the Bessel function from  this note.

In the paper \cite{Ob} we study the space of the $W_k$-valued
functions on  $GL(2n,\CC)$ which are invariant with respect to the
left action of $GL(n,\CC)\times GL(n,\CC)$ and
$W_{k,c}$-equivariant with respect to the right action of
$GL(n,\CC)\times GL(n,\CC)$. In \cite{Ob} we use representation
theory of $GL(2n,\CC)$ to construct the Macdonald-Koornwinder
polynomials. If we degenerate $GL(2n,\mathbb{C})$ into the group
$\mathfrak{gl}(n,\CC)^{\oplus 2}\rtimes GL(n,\CC)^{\times 2}$, we
get  Corollary~\ref{maincor}  for $\ell=2$. Indeed, the group
$\mathfrak{gl}(n,\CC)^{\oplus 2}\rtimes GL(n,\CC)^{\times 2}$ is
the complexification of the group $\bold{H}_n(\mathbb{R})$ and the
right and left action of $GL(n,\CC)^{\times 2}$ on $GL(2n,\CC)$
degenerates into the right and left action of  $GL(n,\CC)^{\times
2}$ on the complexification of the group $\bold{H}_n(\mathbb{R})$.

{\bf Acknowledgements.} I am very grateful to my advisor Pavel
Etingof for the formulation of the problem and many fruitful
discussions which led me to the solution. I  want to thank
Wee Liang Gan for useful discussions and explaining his results with
Victor Ginzburg. I also want to thank
Iain Gordon for the attention to the work and pointing out the
misprints in the preliminary
version of the note. Finally, I am very grateful to Jasper Stokman for
reading the note and giving many useful references.


\begin{thebibliography}{9}
\bibitem[BEG]{BEG} Yu. Berest, P. Etingof, V. Ginzburg, Cherednik algebras
and differential operators on quasi-invariants, Duke Math. J. 118,
no. 2 (2003), 279-337.
\bibitem[BEG1]{BEG1}Yu. Berest, P. Etingof, V. Ginzburg, Finite dimensional
representations of rational Cherednik algebras, IMRN, no. 19 (2003), 1053-1088.
\bibitem[Ch]{Ch} I. Cherednik, Double affine Hecke algebras, KZ equations, and
Macdonald operators, IMRN, no. 9 (1992), 171-180.
\bibitem[DO]{DO} C. Dunkl, E. Opdam, Dunkl operators for
complex reflection groups. Proc. London Math. Soc.
no. 3 86 (2003), no. 1, 70--108.
\bibitem[EFK]{EFK} P.I.Etingof, I.B.Frenkel, A.A.Kirillov, Jr.,
Spherical functions on affine Lie group, Duke Math. J.,  80,  no. 1, (1995),
 59--90.
\bibitem[E]{E} P. Etingof, Cherednik and Hecke algebras of varieties with a finite group action,
math.QA/0406499.
\bibitem[EG]{EG}P. Etingof, V. Ginzburg, Symplectic reflection algebras,
Calogero-Moser space, and deformed Harish-Chandra homomorphism,
math.AG/0011114, Invent. Math., 147 no. 2, (2002), 243-348.
\bibitem[Ga]{Gan} W. L. Gan, Chavalley restrictriction theorem for the cyclic quiver,
math.RT/0405103.
\bibitem[GG]{GG} W. L. Gan, V. Ginzburg, Almost-commuting variety, $\mathcal{D}$-modules,
and Cherednik algebras, math.RT/0409262.

\bibitem[G]{G} I. Gordon,
A remark on rational Cherednik algebras and differential operators on the cyclic quiver,
to appear Glasgow Math Journal.

\bibitem[Ho]{H} M. Holland, Quantization of the Marsden-Weinstein reduction for extended
Dynkin quivers, Ann. Sci. Ecole Norm. Sup., 32 (1999), 813--834.
\bibitem[K]{K} A. A. Kirillov, Lectures on the orbit Method,
Graduate Studies in Mathematics, vol. 64, American Mathematical Society, Providence,
R. I.
\bibitem[LS]{LS}  T. Levasseur, J. Stafford,  Invariant differential operators
and an homomorphism of Harish-Chandra, J. Amer. Math. Soc. 8, no. 2 (1995), 365--372.
\bibitem[Ob]{Ob} A. Oblomkov,  Heckman-Opdam's Jacobi polynomials for the
$BC_n$ root system and generalized spherical functions.
Adv. Math. 186 no. 1, (2004),  153--180.
\bibitem[Op]{Op}E.M. Opdam, Bessel functions and the discriminant of a finite
Coxeter group. Compositio Math. 85, no. 3 (1993),  333--373.
\bibitem[V]{V} N. Ja. Vilenkin, N. Ja. Special functions and the
theory of group representations,
Translations of Mathematical Monographs,
Vol. 22 American Mathematical Society, Providence, R. I.
\end{thebibliography}
\end{document}